\documentclass[12pt,a4paper]{article}
\usepackage{eepic,epic,url,amsfonts, color}
\usepackage{graphicx, graphics}
\usepackage{pstricks,pst-node,pst-tree}
\usepackage{epic}
\usepackage{isorot}
\begin{document}

\newtheorem{theorem}{Theorem}[section]
\newtheorem{prop}[theorem]{Proposition}
\newtheorem{cor}[theorem]{Corollary}
\newtheorem{lemma}[theorem]{Lemma}
\newtheorem{ex}[theorem]{Example}
\newtheorem{no}[theorem]{Note}
\newtheorem{defn}[theorem]{Definition}
\newtheorem{unnumber}{}
\renewcommand{\theunnumber}{\relax}
\newtheorem{prepf}[unnumber]{Proof}
\newenvironment{pf}{\prepf\rm}{\endprepf}
\newcommand{\qed}{\qquad$\color{black}\star$}

\title{General combinatorical structure of truth tables of bracketed formulae connected by implication.}
\author{Volkan Yildiz \\\\
\texttt{vo1kan@hotmail.co.uk}\\
Or\\
\texttt{ah05146@qmul.ac.uk}}
\date{}
\maketitle

\begin{center}
\color{pink}\framebox{\color{pink}\framebox{\parbox[b]{14 cm}{{\includegraphics[scale=0.33]{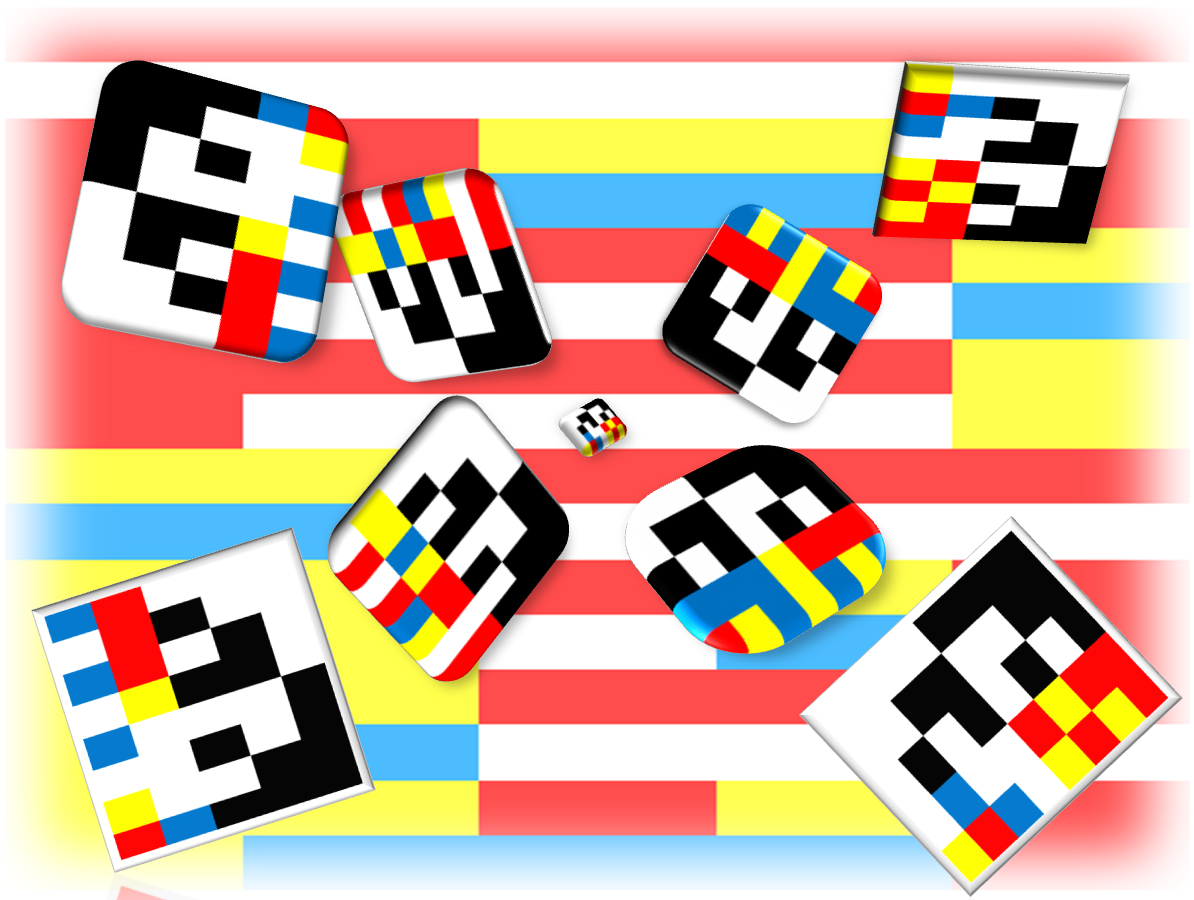}}}}}
\end{center}

\begin{abstract}
In this paper we investigate the general combinatorical structure of the 
truth tables of all bracketed formulae with $n$ distinct variables 
connected by the binary connective of implication, an m-implication.
\end{abstract}

{\footnotesize
{\em Keywords:} 
Propositional logic, implication, Catalan numbers, parity, asymptotic.

AMS classification: 05A15, 05A16, 03B05, 11B75}

\pagebreak

\newcommand
\HRule{\noindent\rule{\linewidth}{0.5pt}}

{\small
\begin{flushleft}
{\em Life is a tragedy when seen in close-up, but a comedy in long-shot. }
\end{flushleft}
\begin{flushright}
{\em C. Chaplin}
\end{flushright}
}
\HRule
$\;$\\
\HRule

\tableofcontents 

\pagebreak

{\bf Notation }\\
{\small
\begin{itemize}
\item[(1)] $p_1,p_2,\ldots , p_n$ and $\phi, \psi$ are propositional variables.
\item[(2)] `True' will be denoted by 1
\item[(3)] `False' will be denoted by 0 
\item[(4)] Set of counting numbers is denoted by $\mathbb{N}=\{ 1,2,3,4,\ldots \}$
\item[(5)] Set of even counting numbers is denoted by $\mathbb{E}$
\item[(6)] Set of odd counting numbers is denoted by $\mathbb{O}$
\item[(7)]$:\quad $ such that
\item[(8)] $\nu$ is the valuation function $:\quad \nu(\phi)=1$ if $\phi$ is true, and $\nu(\phi)=0$ if $\phi$ is false.
\item[(9)] $\wedge$ the conjunction 
\item[(10)] $\vee$ the disjunction
\item[(11)] $ \to $ the binary connective of implication
\item[(12)] $\neg$ the negation operator
\item[(13)] $\# c$ denotes the case number in $t_n^{\#c}$,  where $c=1,2,3$.
\item[(14)] \qed $\;\;\;$ is for QED.
\end{itemize}
}

\HRule
$\;$
\HRule

\section{Preface}

This project begun with the following question: {\em Can we count the number of true entries in truth 
tables of bracketed formulae connected by implication, and modified-implication rules?} 
The primitive answer is to get the number of false entries then subtract it from the total number of entries.
 But this calculation is not precise enough to see the fruitfulness of the truth tables of these kind. 
In this project we have underlined that `unlike the false entries in the truth tables that are connected by 
the binary connective of implication or by m-implication, true entries are not homogeneous structures'. 
The cover picture on the title page is designed to give an intuitive background for this inhomogeneous structure.
\\\\
By reading this project the reader will encounter eight Catalan like sequences, and their asymptotics. 
We have started this project by counting the total number of entries (or rows), in all truth tables for bracketed implication.
Then in section 3, we used the recurrence relation of the total number of rows and we worked backwards to
find out the general structure of the truth tables. Alas, we showed that the structure of the recurrence relation 
of the true, or false entries is inherited from the structure of the recurrence relation of the total number of entries. In section 
4  and  6 we counted the true entries first in using ordinary implication, and later in modified 
binary connective of implication. In section 5 and 7 we dealt with asymptotic of the 
sequences that we have seen in section 4 and 6 respectively. In the last section 
we showed that the sequences that are in section 4 and 6 preserve the parity of Catalan numbers.

\begin{flushright} Volkan Yildiz,\\ 24-May-2012, \\London. \end{flushright}

\pagebreak

\section{Introduction}

In \cite{P1} and \cite{P3} we have mentioned that the total number of rows in all truth tables of 
bracketed implications with $n$ distinct variables is $g_n$ and 
it has the following generating function and the explicit formula respectively:

\begin{center}
\[
G(x)=\frac{1-\sqrt{1-8x}}{2},\quad \textit{ and } \quad g_n=2^nC_n, 
\]
\end{center}

\[
\quad \quad  \textit{ where } \; C_n \textit{ is the nth Catalan number and } \quad  C_n =\frac{1}{n}{{2n-2}\choose{n-1}}.
\]
By using the explicit formula it is straightforward to calculate the values of $g_n$. 
The table below illustrates this up to $n=12$.
\medskip{\small
\[
\begin{array}{|l|c|c|c|c|c|c|c|c|c|c|c|c|c|}
\hline n &0 & 1& 2 & 3 & 4 & 5 & 6 & 7 & 8 & 9 & 10 & 11 & 12 \\
\hline g_n & 0 & 2 & 4 & 16 & 80 & 448 & 2688 & 16896 & 109824 & 732160 & 4978688 & 34398208 & 240787456 \\
\hline
\end{array}
\]
}

\begin{prop}Let $g_n$ be the total number of rows in all truth 
tables for bracketed implication with $n$ distinct variables $p_1, \ldots ,p_n$. Then
\begin{center}
\framebox{\parbox[b]{8 cm}{{\begin{equation}\label{eq;1}
g_n= \sum_{i=1}^{n-1} g_ig_{n-i}, \;\;\;\;\; \textit{ with } \;\;\; g_1=2.
\end{equation}}}

}
\end{center}
\end{prop}
 
\begin{pf}
[First Proof]\\\
For $p_1\to \ldots \to p_n$ there are $2^n$ rows and $C_n$ columns, hence there are $2^nC_n$ rows altogether. 
\[
g_n= 2^nC_n= 2^n\sum_{i=1}^{n-1}C_iC_{n-i} = \sum_{i=1}^{n-1}(2^iC_i)(2^{n-i}C_{n-i}) 
= \sum_{i=1}^{n-i} g_ig_{n-i}, \; \textit{ with } g_1=2. \;\; \star
\]
[Second proof, `intuitive']  \\\\
Consider $n$ distinct propositions $p_1, \ldots , p_n$, for $i\geq 1$:

\[
\underbrace{p_1\to \ldots \to p_i}_{g_i \textit{rows}}\to \underbrace{p_{i+1}\to \ldots \to p_n}_{g_{n-i} \textit{rows}} 
\]
 There are $g_i$ rows for $p_1\to \ldots \to p_i$, and there are $g_{n-i}$ rows for $p_{i+1}\to \ldots \to p_n$, 
and the number of choices is given the recurrence relation~(\ref{eq;1}) .\qed
\end{pf}

\section{General Structure}
To see the combinatorical structure of the truth tables of bracketed formulae connected by implication, we need to 
make use of the method of working backwards. This method allows us to avoid unnecessary choices altogether.
In ancient Greek it was mentioned by the mathematician Pappus and in recent times the method of working backwards has
been discussed G. Polya. 
\medskip
Recall that
\[
g_n=\sum_{i=1}^{n-1}g_ig_{n-i}
\]
Since each table consists of false and true rows, we let $g_i=t_i+f_i$, 
where $t_i$, $f_i$ is the corresponding number of `true', `false' rows in $g_i$ respectively. Then

\[
g_n= \sum_{i=1}^{n-1} (t_i+f_i)(t_{n-i}+f_{n-i}) 
\]
If we expand the right hand side:
\begin{equation} \label{eq;2}
g_n= \underbrace{\sum_{i=1}^{n-1}t_i t_{n-i}}_{\textit{case 1}} + \underbrace{\sum_{i=1}^{n-1}f_i t_{n-i}}_{\textit{case 2}} + \underbrace{\sum_{i=1}^{n-1}f_i f_{n-i}}_{\textit{case 3}} + \underbrace{\sum_{i=1}^{n-1}t_i f_{n-i}}_{\textit{case 4}}  
\end{equation}

We can now partition the right hand side into more tangible cases and explain what each of these summands mean. 
Let $\psi$ and $\phi$ be propositional variables, and let $\nu$ be the {\em valuation} function, 
then 
\[
\nu(\psi \to \phi) = 1 \Longleftrightarrow (\nu(\psi)=0 \vee \nu(\phi)=1).
\]
Therefore there are three cases to consider here:

\[
 \underbrace{(\nu(\psi)=1 \wedge \nu(\phi)=1)}_{\textit{ case 1}} \vee \underbrace{(\nu(\psi)=0 \wedge \nu(\phi)=1)}_{\textit{case 2}} 
\vee \underbrace{(\nu(\psi)=0 \wedge \nu(\phi)=0)}_{\textit{case 3}}
\]
Addition to the three cases in above there is also the fourth case, known as the {\em `disastrous combination'} :
\[
\nu(\psi \to \phi) = 0 \Longleftrightarrow \underbrace{(\nu(\psi)=1 \wedge \nu(\phi)=0)}_{\textit{case 4}}.
\]
\medskip
The four cases in equation~(\ref{eq;2}) coincide with the four cases in the penultimate lines respectively, 
and we get the following theorem:
\begin{theorem} \label{T:1}
Let $g_n$ be the total number of rows in all truth tables of 
bracketed implications with $n$ distinct variables. Then for $n\geq 2$, $g_n$ can be partition 
into four cases as below
\[
g_n= t_n^{\#1}+t_n^{\#2}+t_n^{\#3}+f_n, \textit{ where }
\]
\[
 t_n^{\#1}=\sum_{i=1}^{n-1}t_i t_{n-i}, \quad  t_n^{\#2}=\sum_{i=1}^{n-1}f_i t_{n-i},\quad  t_n^{\#3}=\sum_{i=1}^{n-1}f_i f_{n-i}, \textit{ and} \quad  f_n=\sum_{i=1}^{n-1}t_i f_{n-i}. 
\]
with
\[
0=t_1^{\#1}=t_1^{\#2}=t_1^{\#3}, \quad \textit{ and }\quad  f_1=1.
\]
\end{theorem}
In coming parts of this paper we will investigate each of the above cases in more detail.

\section{\small Counting true entries in truth tables of
bracketed formulae connected by implication}

\subsection{Case 4}
Let $\psi$ and $\phi$ be propositional variables, then

\[
\nu(\psi\to \phi)=0 \quad : \quad (\nu(\psi )= 1 \wedge \nu(\phi )=0).
\]
This case has been discussed already in  \cite{P1}; we had the following results.
\begin{theorem}
Let $f_n$ be the number of rows with the value ``false'' in the truth tables 
of all bracketed formulae with $n$ distinct propositions $p_1,\ldots,p_n$ 
connected by the binary connective of implication. Then 

\begin{center}
\framebox{\parbox[b]{8 cm}{{\begin{equation}\label{eq1:f}
f_n =\sum_{i=1}^{n-1} (2^iC_i-f_i)f_{n-i}, \; \textit{ with } \; f_1=1
\end{equation}}}}
\end{center}

and for large $n$, $\;f_n \sim \Bigg(\frac{3-\sqrt{3}}{6}\Bigg)\frac{2^{3n-2}}{\sqrt{\pi n^3}}$. 
Where $C_i$ is the $i$th Catalan number.
\end{theorem}
If we look closely to the recurrence relation~(\ref{eq1:f}), since $t_i= (2^iC_i-f_i)$ 
we obtain the `case 4' in equation~(\ref{eq1:f}):
\[
f_n =\sum_{i=1}^{n-1} (2^iC_i-f_i)f_{n-i}=\sum_{i=1}^{n-1} t_if_{n-i}, \; \textit{ with } \;  f_1=1 .
\]
First ten terms of $\{f_n\}_{n> 0}$ are,
\[
1, 1, 4, 19, 104, 614, 3816, 24595, 162896, 1101922,\ldots
\]
We have also shown in \cite{P1} that $f_n$ has the following generating function:
\begin{center}
\framebox{\parbox[b]{10 cm}{{
\[
F(x) = \frac{-1-\sqrt{1-8x} + \sqrt{2+2\sqrt{1-8x}+8x}}{4}.
\]
}}}
\end{center}
We have also shown in \cite{P2} that the sequence $\{f_n\}_{n\geq 1}$ preserves the parity of Catalan numbers.

\subsection{Case 3}

Let $\psi$ and $\phi$ be propositional variables, then
\[\nu(\psi\to \phi)=1 \quad  : \quad (\nu(\psi) = 0 = \nu(\phi)) .\]

In this case we are interested in formulae obtained from $p_1\to \ldots \to p_n$ by inserting brackets so that 
the valuation of the first $i$ bracketing and the rest $(n-i)$ bracketing both give 0, `false'. 
The table 1 below shows the truth tables, (merged into one), for the two bracketed implications
in $n=3$ variables. The corresponding case 3 truth values are denoted in green.
{\small
\begin{table}[h]
\caption{$n=3$}
\centering
\begin{tabular}{|c| c | c | c | c|}
\hline \hline
$p_1$ & $p_2$ & $p_3$ & $p_1\to (p_2 \to p_3)$ & $(p_1 \to p_2)\to p_3$ \\
\hline
1 & 1 & 1 & 1 & 1 \\
\hline 
1 & 1 & 0 &  0 &  0\\
\hline 
1 & 0 & 1 & 1 & 1 \\
\hline 
1 & 0 & 0 & 1 & {\bf \color{green} 1} \\
\hline 
0 & 1 & 1 & 1 & 1 \\
\hline 
0 & 1 & 0 & {\bf \color{green} 1} & 0 \\
\hline 
0 & 0 & 1 & 1 & 1 \\
\hline 
0 & 0 & 0 & 1 & 0\\
\hline
\end{tabular}
\end{table} 
}

\begin{prop}
Let $t_n^{\#3}$ be the number of rows with the value ``true'' in the truth tables 
of all bracketed formulae with $n$ distinct propositions $p_1,\ldots,p_n$ 
connected by the binary connective of implication such that the valuation 
of the first $i$ bracketing and the rest $(n-i)$ bracketing both  give 0, `false'. Then
\begin{equation}\label{e:t3}
t_n^{\#3} =\sum_{i=1}^{n-1} f_if_{n-i} \; \textit{ where } \{f_i\}_{i\geq 1} \textit{ is the sequence in equation~(\ref{eq1:f}). }  
\end{equation}
\end{prop}
\begin{pf}
A row with the value 1, `true', comes from an expression $\psi\to \phi$, where $\nu(\psi)=0$ and $\nu(\phi)=0$. 
If $\psi$ contains $i$ variables, then $\phi$ contains $(n-i)$ variables, and the number of 
choices is given by the summand in the proposition. \qed
\end{pf}

It is now very easy to find out the generating function of $t_n^{\#3}$ by using the relation~(\ref{e:t3}).
 Let $T_3(x)=\sum_{n\geq 1} t_n^{\#3} x^n$ then  
$T_3(x) = F(x)^2$, and hence we get the following proposition:
\begin{prop}
The generating function for the sequence $\{ t_n^{\#3} \}_{n\geq 0}$ is given by 

\begin{center}
\framebox{\parbox[b]{16 cm}{{
\[
T_3(x)= \frac{2+2\sqrt{1-8x}-\sqrt{2+2\sqrt{1-8x}+8x}-\sqrt{1-8x}\sqrt{2+2\sqrt{1-8x}+8x}}{8} .
\]
}}}
\end{center}

\end{prop}

By using Maple we find the first 21 terms of this sequence:

\begin{eqnarray*}
\{t_n^{\#3}\}_{n>1} &=&  1, 2, 9, 46, 262, 1588, 10053, 65686, 439658, 2999116, \\
&& 20774154, 145726348, 1033125004, 7390626280, 53281906861, \\
&& 386732675046, 2823690230850, 20725376703324, \\
&& 152833785130398, 1131770853856100, 8412813651862868, \ldots
\end{eqnarray*}

We will discuss the asymptotic and number theoretical results in a separate chapter. 
Now we have to embark on case 2.

\subsection{Case 2}
Let $\psi$ and $\phi$ be propositional variables, then
\[
\nu(\psi\to \phi)=1 \quad : \quad (\nu(\psi)=0 \; \wedge \; \nu(\phi)=1) .
\]

In this case we are interested in formulae obtained from $p_1\to \ldots \to p_n$ by inserting brackets so that 
the valuation of the first $i$ bracketing give 0, `false', and the rest $(n-i)$ bracketing  give 1, `true'.

\begin{prop}
Let $t_n^{\#2}$ be the number of rows with the value ``true'' in the truth tables 
of all bracketed formulae with $n$ distinct propositions $p_1,\ldots,p_n$ 
connected by the binary connective of implication such that the valuation 
of the first $i$ bracketing is 0 and and the rest $(n-i)$ bracketing gives 1. Then,
\begin{equation}\label{e:t2}
t_n^{\#2}=\sum_{i=1}^{n-1}f_it_{n-i}
\end{equation}
\end{prop}
\begin{pf}
A row with the value 1, `true', comes from an expression $\psi\to \phi$, where $\nu(\psi)=0$ and $\nu(\phi)=1$. 
If $\psi$ contains $i$ variables, then $\phi$ contains $(n-i)$ variables, and the number of 
choices is given by the summand in the proposition. \qed
\end{pf}

The table 2 above, (merged into one), for the two bracketed implications
in $n=3$ variables, indicates the corresponding truth values in bold.

{\small
\begin{table}[h]
\caption{$n=3$}
\centering
\begin{tabular}{|c| c | c | c | c|}
\hline \hline
$p_1$ & $p_2$ & $p_3$ & $p_1\to (p_2 \to p_3)$ & $(p_1 \to p_2)\to p_3$ \\
\hline
1 & 1 & 1 & 1 & 1 \\
\hline 
1 & 1 & 0 &  0 &  0\\
\hline 
1 & 0 & 1 & 1 & {\bf 1} \\
\hline 
1 & 0 & 0 & 1 & 1 \\
\hline 
0 & 1 & 1 & {\bf 1} & 1 \\
\hline 
0 & 1 & 0 & 1 & 0 \\
\hline 
0 & 0 & 1 & {\bf 1} & 1 \\
\hline 
0 & 0 & 0 & {\bf 1} & 0\\
\hline
\end{tabular}
\end{table} 
}

\begin{cor}
$t_n^{\#2}= f_n$ for $n\geq 2$ .
\end{cor}
\begin{pf}
Since $t_{n-i}=(2^{n-i}C_{n-i}-f_{n-i})$,
\begin{equation}\label{e:t22}
t_n^{\#2}=\sum_{i=1}^{n-1}f_i(2^{n-i}C_{n-i}-f_{n-i})
\end{equation}
More explicitly,
\begin{eqnarray*}
t_n^{\#2}&=&  f_1(2^{n-1}C_{n-1}-f_{n-1})+\ldots + f_{n-2}(2^{2}C_{2}-f_{2})+f_{n-1}(2^{1}C_{1}-f_{1}) \\
&=& f_{n-1}(2^{1}C_{1}-f_{1})  +  f_{n-2}(2^{2}C_{2}-f_{2})+\ldots +  f_1(2^{n-1}C_{n-1}-f_{n-1})\\
&=& \sum_{i=1}^{n-1} (2^iC_i-f_i)f_{n-i} \\
&=& f_n.
\end{eqnarray*}
\qed
\end{pf}
Thus  for $n\geq 2$,  Case 2 coincides with Case 4, which we have investigated in great detail in \cite{P1}, and in \cite{P2}. 

\begin{prop}
The generating function for the sequence $\{ t_n^{\#2} \}_{n\geq 0}$ is given by 

\begin{center}
\framebox{\parbox[b]{14 cm}{{\[
T_2(x)= F(x)-x = \frac{-1-\sqrt{1-8x} + \sqrt{2+2\sqrt{1-8x}+8x} -4x}{4} .
\]}}}
\end{center}

\end{prop}

\subsection{Case 1}

Let $\psi$ and $\phi$ be propositional variables, then
\[
\nu(\psi\to \phi)=1 \quad : \quad (\nu(\psi)=1= \nu(\phi)) .
\]

In this case we are interested in formulae obtained from $p_1\to \ldots \to p_n$ by inserting brackets so that 
the valuation of the first $i$ bracketing  and the rest $(n-i)$ bracketing  give 1, `true'. 

Table 3, (merged into one), for the two bracketed implications
in $n=3$ variables, indicates the corresponding truth values in red.

{\small
\begin{table}[h]
\caption{$n=3$}
\centering
\begin{tabular}{|c| c | c | c | c|}
\hline \hline
$p_1$ & $p_2$ & $p_3$ & $p_1\to (p_2 \to p_3)$ & $(p_1 \to p_2)\to p_3$ \\
\hline
1 & 1 & 1 & {\bf\color{red}1} & {\bf\color{red}1} \\
\hline 
1 & 1 & 0 &  0 &  0\\
\hline 
1 & 0 & 1 & {\bf\color{red}1} &  1 \\
\hline 
1 & 0 & 0 & {\bf\color{red}1} & 1 \\
\hline 
0 & 1 & 1 & 1 & {\bf\color{red}1} \\
\hline 
0 & 1 & 0 & 1 & 0 \\
\hline 
0 & 0 & 1 & 1 & {\bf\color{red}1} \\
\hline 
0 & 0 & 0 & 1 & 0\\
\hline
\end{tabular}
\end{table} 
}

\begin{prop}
Let $t_n^{\#1}$ be the number of rows with the value ``true'' in the truth tables 
of all bracketed formulae with $n$ distinct propositions $p_1,\ldots,p_n$ 
connected by the binary connective of implication such that the valuation 
of the first $i$ bracketing  and and the rest of $(n-i)$ bracketing gives 1. 
Then,
\begin{equation}\label{eq:t2}
t_n^{\#1}=\sum_{i=1}^{n-1}t_it_{n-i}
\end{equation}

\end{prop}
\begin{pf}
A row with the value 1, `true', comes from an expression $\psi\to \phi$, where $\nu(\psi)=1$ and $\nu(\phi)=1$. 
If $\psi$ contains $i$ variables, then $\phi$ contains $(n-i)$ variables, and the number of 
choices is given by the summand in the proposition. \qed
\end{pf}

Let $T_1(x)=\sum_{n\geq 1}t_n^{\#1} x^n$
, since $G(x)=F(x)+T_1(x)+T_2(x)+T_3(x)\; $ then 
it is easy to find that $T_1(x)$ has the following explicit algebraic representation:

\begin{prop} The generating function for the sequence $\{ t_n^{\#1} \}_{n\geq 0}$ is given by

\begin{center}
\framebox{\parbox[b]{16 cm}{{\[
T_1(x)=\frac{6-2\sqrt{1-8x}-3\sqrt{2+2\sqrt{1-8x}+8x}+ \sqrt{1-8x}\sqrt{2+2\sqrt{1-8x}+8x}}{8} .
\]}}}
\end{center}

\end{prop}

\medskip

By using Maple we find the first 21 terms of this sequence:
\medskip

\begin{eqnarray*}
\{t_n^{\#1} \}_{n\geq 2} &=& 1, 6, 33, 194, 1198, 7676, 50581, 340682, 2335186, 16237284, \\
&& 114255994, 812107412, 5822171548, 42052209400, 305714145869, \\
&& 2235262899418, 16426616425002, 121265916776148, \\
&& 898878250833358, 6687497426512700, 49920590244564484, \ldots
\end{eqnarray*}

The below table shows the sequences which we have discussed so far, up to $n=11$.

{\small
\medskip{\small
\[
\begin{array}{|l|c|c|c|c|c|c|c|c|c|c|c|c|}
\hline n &0 & 1& 2 & 3 & 4 & 5 & 6 & 7 & 8 & 9 & 10 & 11 \\
\hline 2^n & 1 & 2 & 4 & 8 & 16 & 32 & 64 & 128 & 256 & 512 &  1024 & 2048  \\
\hline C_n & 0 & 1 & 1 & 2 & 5 & 14 & 42 & 132 & 429 & 1430 & 4862 & 16796  \\
\hline g_n & 0 & 2 & 4 & 16 & 80 & 448 & 2688 & 16896 & 109824 & 732160 & 4978688 & 34398208  \\
\hline f_n & 0 & 1 & 1 & 4 & 19 & 104 & 614 & 3816 & 24595 & 162896 & 1101922 & 7580904  \\
\hline t_n^{\#1} & 0 & 0 & 1 & 6 & 33 & 194 & 1198 & 7676 & 50581 & 340682 & 2335186 & 16237284 \\
\hline t_n^{\#2} & 0 & 0 & 1 & 4 & 19 & 104 & 614 & 3816 & 24595 & 162896 & 1101922 & 7580904  \\
\hline t_n^{\# 3} & 0 & 0 & 1 & 2 & 9 & 46 & 262 & 1588 & 10053 & 65686 & 439658 & 2999116 \\
\hline
\end{array}
\]
}

\begin{sidewaystable}

\caption{ The below truth tables, (merged into one), for the five bracketed implications
in $n=4$ variables. Where Case 1 is in red, case 2 in black, case 3 in green and case 4 is indicated in blue.}
\centering
\begin{tabular}{|c| c | c | c | c|c|c|c|c|}
\hline \hline
$p_1$ & $p_2$ & $p_3$ & $p_4$ & $p_1\to (p_2\to (p_3 \to p_4))$ & $p_1 \to(( p_2\to p_3)\to p_4)$ & $(p_1 \to (p_2\to p_3))\to p_4$ & $ ((p_1 \to p_2)\to p_3)\to p_4$ &  $ (p_1 \to p_2)\to (p_3\to p_4)$\\
\hline
1&1&1&1 & {\bf \color{red} 1} & {\bf \color{red} 1} & {\bf \color{red} 1} & {\bf \color{red} 1} & {\bf \color{red} 1} \\
\hline
1&1&1&0 & {\bf \color{blue} 0} & {\bf \color{blue} 0}& {\bf \color{blue} 0}& {\bf \color{blue} 0}& {\bf \color{blue} 0}\\
\hline
1&1&0&1& {\bf \color{red} 1} & {\bf \color{red} 1} & {\bf  1} & {\bf  1} & {\bf \color{red} 1} \\
\hline
1&1&0&0& {\bf \color{red} 1} & {\bf \color{red} 1} & {\bf \color{green} 1} & {\bf \color{green} 1}&  {\bf \color{red} 1} \\
\hline
1&0&1&1& {\bf \color{red} 1} & {\bf \color{red} 1} & {\bf \color{red} 1} & {\bf \color{red} 1} & {\bf  1} \\
\hline 
1&0&1&0& {\bf \color{red} 1} & {\bf \color{blue} 0}& {\bf \color{blue} 0}& {\bf \color{blue} 0} & {\bf \color{green} 1}\\
\hline
1&0&0&1& {\bf \color{red} 1} & {\bf \color{red} 1} & {\bf \color{red} 1} & {\bf \color{red} 1} & {\bf  1} \\
\hline 
1&0&0&0& {\bf \color{red} 1} & {\bf \color{blue} 0} & {\bf \color{blue} 0} & {\bf \color{blue} 0} & {\bf 1} \\
\hline
0&1&1&1 & {\bf 1} & {\bf 1} & {\bf \color{red} 1} & {\bf \color{red} 1} & {\bf \color{red} 1} \\
\hline
0&1&1&0 & {\bf\color{green} 1} & {\bf\color{green} 1} & {\bf \color{blue} 0} & {\bf \color{blue} 0}& {\bf \color{blue} 0}\\
\hline
0&1&0&1 & {\bf 1} & {\bf 1} & {\bf\color{red} 1} & {\bf 1} & {\bf \color{red} 1} \\
\hline
0&1&0&0 & {\bf 1} & {\bf 1} & {\bf\color{blue} 0} & {\bf \color{green} 1} & {\bf \color{red} 1} \\
\hline
0&0&1&1 & {\bf 1} & {\bf 1} & {\bf\color{red} 1} & {\bf \color{red} 1} & {\bf \color{red} 1} \\
\hline 
0&0&1&0 & {\bf 1} & {\bf \color{green} 1} & {\bf\color{blue} 0} & {\bf\color{blue} 0} & {\bf\color{blue} 0} \\
\hline
0&0&0&1 & {\bf 1} & {\bf 1} & {\bf \color{red} 1}& {\bf 1}&{\bf \color{red} 1}\\
\hline
0&0&0&0 & {\bf 1} & {\bf\color{green} 1} & {\bf\color{blue} 1} &  {\bf\color{green} 1}  & {\bf \color{red} 1}\\
\hline

\end{tabular}
\label{e:t}
\end{sidewaystable}

\section{Asymptotic Estimate 1}
In this chapter we will mainly make use of the asymptotic techniques that we have 
mentioned in~\cite[pg 6-7]{P3}, and in~\cite[pg 6-7]{P1}.  
In \cite{P1} we have shown that the following asymptotic results 
for the sequence $\{f_n\}_{n\geq 1}$ are true:

\begin{theorem}
Let $f_n$ be number of rows with the value false in the truth tables of all the
bracketed implications with $n$ variables. Then

\begin{center}
\framebox{\parbox[b]{6 cm}{{\[
f_n \sim \left(\frac{3-\sqrt{3}}{6}\right)\frac{2^{3n-2}}{\sqrt{\pi n^3}}.
\]
}}}
\end{center}

\end{theorem}

\begin{cor}
Let $g_n$ be the total number of rows in all truth tables for bracketed
implications with $n$ variables, and $f_n$ the number of rows with the
value ``false''. Then $\lim_{n\to\infty}f_n/g_n=(3-\sqrt{3})/6$.
\end{cor}
 
\begin{theorem}
Let $t_n^{\#3}$ be number of rows with the value true in the truth tables of all the
bracketed implications with $n$ variables, case 3. Then

\begin{center}
\framebox{\parbox[b]{6 cm}{{\[
t_n^{\#3} \sim \left(\frac{2\sqrt{3}-3}{6}\right)
\frac{2^{3n-2}}{\sqrt{\pi n^3}}
\]}}}
\end{center}
\end{theorem}

\begin{pf}
Recall that 

{\small
\[
T_3(x)= \frac{2+2\sqrt{1-8x}-\sqrt{2+2\sqrt{1-8x}+8x}-\sqrt{1-8x}\sqrt{2+2\sqrt{1-8x}+8x}}{8} .
\]
}
By using the asymptotic techniques which we have discussed in \cite{P1} and \cite{P3}, we found that
$r=\frac{1}{8}$, $f(x)= \sqrt{1-8x}$ and $T_3(\frac{1}{8})= \frac{2-\sqrt{3}}{8}\not= 0$. Let
$A(x)=T_3(x)-T_3(\frac{1}{8})$

\[\lim_{x\to1/8} \frac{A(x)}{f(x)} 
= \lim_{x\to 1/8} \frac{2\sqrt{1-8x}-\sqrt{2+2\sqrt{1-8x}+8x}-\sqrt{1-8x}\sqrt{2+2\sqrt{1-8x}+8x}+\sqrt{3}}{8\sqrt{1-8x}} .\]

Let $v=\sqrt{1-8x}$. Then
\begin{eqnarray*}
L &=& \lim_{v\to 0} \frac{2v -\sqrt{(1+v)(3-v)}-v\sqrt{(1+v)(3-v)} +\sqrt{3}}{8v} \\
&=& \lim_{v\to 0} \frac{2v -\sqrt{3+2v-v^2} -v\sqrt{3+2v-v^2}+ \sqrt{3}}{8v} \\
&=& \lim_{v\to 0} \frac{2(\sqrt{-(v+1)(v-3)}-2-v+v^2)}{8\sqrt{-(v+1)(v-3)}}\\
&=& - \frac{2\sqrt{3}-3}{12},
\end{eqnarray*}
where we have used l'H\^opital's Rule in the penultimate line.

Finally,
\[t_n^{\#3} \sim - \frac{2\sqrt{3}-3}{12}{n-\frac{3}{2}\choose n}
\left(\frac{1}{8}\right)^{-n} \sim \left(\frac{2\sqrt{3}-3}{6}\right)
\frac{2^{3n-2}}{\sqrt{\pi n^3}},\]
and the proof is finished. \qed
\end{pf}

Recall that
\[
T_2(x)= F(x)-x = \frac{-1-\sqrt{1-8x} + \sqrt{2+2\sqrt{1-8x}+8x} -4x}{4} .
\]

\begin{theorem}
Let $t_n^{\#2}$ be number of rows with the value true in the truth tables of all the
bracketed implications with $n$ variables, case 2. Then

\begin{center}
\framebox{\parbox[b]{6 cm}{{\[
t_n^{\#2} \sim \left(\frac{3-\sqrt{3}}{6}\right)\frac{2^{3n-2}}{\sqrt{\pi n^3}}.
\]}}}
\end{center}

\end{theorem}
\begin{pf} Same as for $f_n$,   see \cite{P1} .\qed \end{pf}

\begin{theorem}
Let $t_n^{\#1}$ be number of rows with the value true in the truth tables of all the
bracketed implications with $n$ variables, case 1. Then

\begin{center}
\framebox{\parbox[b]{6 cm}{{\[
 t_n^{\#1} \sim \left(\frac{1}{2}\right)
\frac{2^{3n-2}}{\sqrt{\pi n^3}}
\]}}}
\end{center}
\end{theorem}

\begin{pf}
Recall that 
\[
T_1(x)=\frac{6-2\sqrt{1-8x}-3\sqrt{2+2\sqrt{1-8x}+8x}+ \sqrt{1-8x}\sqrt{2+2\sqrt{1-8x}+8x}}{8}
\]
By using the asymptotic techniques which we have discussed in \cite{P1} and \cite{P3}, we found that
$r=\frac{1}{8}$, $f(x)= \sqrt{1-8x}$ and $T_1(\frac{1}{8})= \frac{6-3\sqrt{3}}{8}\not= 0$. Let
$A(x)=T_1(x)-T_1(\frac{1}{8})$, then

\[ 
\lim_{x\to1/8} \frac{A(x)}{f(x)} 
= \lim_{x\to 1/8} \frac{-2\sqrt{1-8x}-3\sqrt{2+2\sqrt{1-8x}+8x}+ \sqrt{1-8x}\sqrt{2+2\sqrt{1-8x}+8x}+3\sqrt{3}}{8\sqrt{1-8x}} .
\]

Let $v=\sqrt{1-8x}$. Then
\begin{eqnarray*}
L &=& \lim_{v\to 0} \frac{-2v -3\sqrt{(1+v)(3-v)}+v\sqrt{(1+v)(3-v)} +3\sqrt{3}}{8v} \\
&=& \lim_{v\to 0} \frac{-2v -3\sqrt{3+2v-v^2} +v\sqrt{3+2v-v^2}+ \sqrt{3}}{8v} \\
&=& \lim_{v\to 0} -\frac{2(\sqrt{-(v+1)(v-3)}-3v+v^2)}{8\sqrt{-(v+1)(v-3)}}\\
&=& - \frac{1}{4},
\end{eqnarray*}
where we have used l'H\^opital's Rule in the penultimate line.
Finally,
\[t_n^{\#1} \sim - \frac{1}{4}{n-\frac{3}{2}\choose n}
\left(\frac{1}{8}\right)^{-n} \sim \left(\frac{1}{2}\right)
\frac{2^{3n-2}}{\sqrt{\pi n^3}},\]
and the proof is finished. \qed
\end{pf}

The importance of the constants $\frac{2\sqrt{3}-3}{6}= 0.077350269189$ ,  $\frac{3-\sqrt{3}}{6} = 0.211324865404$ and $\frac{1}{2}=0.5$ lies in the following fact:

\begin{cor}
Let $g_n$ be the total number of rows in all truth tables for bracketed
implications with $n$ variables, then 

\[
\lim_{n\to\infty}\frac{t_n^{\# 3}}{g_n}=\frac{2\sqrt{3}-3}{6}, \;\;\;  \lim_{n\to\infty}\frac{f_n}{g_n}=\lim_{n\to\infty}\frac{t_n^{\# 2}}{g_n}=\frac{3-\sqrt{3}}{6}, \;\;\; \lim_{n\to\infty}\frac{t_n^{\# 1}}{g_n}=\frac{1}{2}. 
\]
\end{cor}
The table below illustrates the convergence :

\[
\begin{array}{|c|c|c|c|c|c|c|c|}
\hline n & f_n = t_n^{\#2} & g_n & t_n^{\#2}/g_n & t_n^{\#1} & t_n^{\#1}/g_n & t_n^{\#3} & t_n^{\#3}/g_n  \\
\hline 1 & 1 \;\; \vline \;\; -  & 2 & 0.5 \;\; \vline \;\;-  & - & - & - & -\\
\hline 2 & 2 & 4 & 0.25& 1& 0.25 & 1 & 0.25 \\
\hline 3 & 4 & 16 & 0.25& 6 & 0.375 & 2 & 0.125\\
\hline 4 & 19 & 80 & 0.2375 & 33 & 0.4125 & 9 & 0.1125\\
\hline 5 & 104 & 428 & 0.2321428571 & 194 & 0.433035714 & 46 & 0.102678571 \\
\hline 6 & 614 & 2688 & 0.228422619 & 1198 & 0.445684524 & 262 & 0.0974702381 \\
\hline 7 & 3816 & 16896 & 0.2258522727 & 7676 & 0.454308712 & 1588 & 0.0939867424 \\
\hline 8 & 424595 & 109824 & 0.2239492279 & 50581 & 0.460564175 & 10053 & 0.0915373689 \\
\hline 9 & 162896 & 732160 & 0.2224868881 & 340682 & 0.465310861 & 65686 & 0.0897153628 \\
\hline 10 & 1101922 & 4978688 & 0.2213277876 & 2335186 &  0.469036421 & 439658 & 0.0883080040 \\
\hline 100 & - & - & 0.212290865  & - &  0.497093847 & - & 0.0783244229  \\
\hline
\end{array}\]

\begin{cor}
Let $\mathcal{P}^{\#1}=\frac{t_n^{\#1}}{g_n}$, $\mathcal{P}^{\#2}=\frac{t_n^{\#2}}{g_n}$, $\mathcal{P}^{\#3}=\frac{t_n^{\#3}}{g_n}$, and $\mathcal{P}^{\#4}=\frac{f_n}{g_n}$, then for $n>2$,

\[
\mathcal{P}^{\#1} > (\mathcal{P}^{\#2} = \mathcal{P}^{\#4}) >\mathcal{P}^{\#3} .
\]

\end{cor}

\begin{cor}
\[
\lim_{n\to\infty}\frac{t_n^{\# 3}}{t_n^{\#2}}=\frac{\sqrt{3}-1}{2}, \;\;\;  \lim_{n\to\infty}\frac{t_n^{\#2}}{t_n^{\#3}}=1+\sqrt{3}, \;\;\; 
\lim_{n\to\infty}\frac{t_n^{\# 3}}{t_n^{\#1}}=\frac{2\sqrt{3}-3}{3}, \;\;\;
\]
\[
  \lim_{n\to\infty}\frac{t_n^{\# 1}}{t_n^{\#3}}=3+2\sqrt{3}, \;\;\;  
\lim_{n\to\infty}\frac{t_n^{\# 2}}{t_n^{\#1}}=\frac{3-\sqrt{3}}{3}, \;\;\;  \lim_{n\to\infty}\frac{t_n^{\#1}}{t_n^{\#2}}=\frac{3+\sqrt{3}}{2}.  
\]	
\end{cor}

\section{\small Counting true entries in truth tables of
bracketed formulae connected by modified implication}
A number of new enumerative problems arise if we modify the 
binary connective of implication as in below.

\subsection{Type 1}

\begin{defn}
Let $\psi$, and $\phi$ be propositional variables then 
\[
\psi \rightharpoonup \phi \equiv \psi \to \neg \phi .
\]
Thus for any valuation $\nu$,
\[\nu(\psi\rightharpoonup\phi)=\cases{0 & if $\nu(\psi)=1$ and $\nu(\phi)=1$,\cr
1 & otherwise;\cr}\] 
\end{defn}
\[\nu(\psi  \rightharpoonup  \phi)=1 \quad : \quad (\nu(\psi)=0 \vee \nu(\phi)=0) \] 
Therefore there are three cases to consider here:
\[
 \underbrace{(\nu(\psi)=0 \wedge \nu(\phi)=0)}_{\textit{ case 1}} \vee \underbrace{(\nu(\psi)=0 \wedge \nu(\phi)=1)}_{\textit{case 2}} 
\vee \underbrace{(\nu(\psi)=1 \wedge \nu(\phi)=0)}_{\textit{case 3}}
\]
Addition to the three cases in above there is also the fourth case, known as the {\em disastrous combination} :
\[
\nu(\psi \rightharpoonup \phi) = 0 \Longleftrightarrow \underbrace{(\nu(\psi)=1 \wedge \nu(\phi)=1)}_{\textit{case 4}}.
\]
\medskip
Let $y_n$, and $d_n$ be the number of rows with the value ``false'' and ``true'' in the truth tables 
of all bracketed formulae with $n$ distinct propositions $p_1,\ldots,p_n$  
connected by the binary connective of m-implication, type 1, respectively. Then using theorem \ref{T:1}, 
\[
g_n = \sum_{i=1}^{n-1} (y_i+d_i)(y_{n-i}+d_{n-i})
\]
Thus if we expand the right hand side 
\begin{equation} \label{eqg:t2}
g_n=\underbrace{\sum_{i=1}^{n-1}y_i y_{n-i}}_{\textit{case 1}} + \underbrace{\sum_{i=1}^{n-1}y_i d_{n-i}}_{\textit{case 2}}  + \underbrace{\sum_{i=1}^{n-1}d_i y_{n-i}}_{\textit{case 3}}   + \underbrace{\sum_{i=1}^{n-1}d_i d_{n-i}}_{\textit{case 4}}  
\end{equation}
The four cases in equation~(\ref{eqg:t2}) coincide with the four cases in the penultimate lines respectively. 
Let
\[
d_n^{\#1}= \sum_{i=1}^{n-1}y_iy_{n-i}, \quad d_n^{\#2}= \sum_{i=1}^{n-1}y_id_{n-i}, \quad d_n^{\#3}= \sum_{i=1}^{n-1}d_iy_{n-i},\quad y_n= \sum_{i=1}^{n-1}d_id_{n-i} .   
\]
Where
\[
0=d_1^{\#1}=d_1^{\#2}=d_1^{\#3}, \quad \textit{ and }\quad y_1=1 .
\]

\subsubsection{Case 4}

This case has been manifested in \cite{P3}, and in summary we had the following results:
\begin{theorem}\label{t:y}
Let $y_n$ be the number of rows with the value ``false'' in the truth tables 
of all bracketed formulae with $n$ distinct propositions $p_1,\ldots,p_n$ 
connected by the binary connective of m-implication, type 1, case 4. Then 
\begin{center}
\framebox{\parbox[b]{12 cm}{{\begin{equation}\label{e:y}
y_n=\sum_{i=1}^{n-1} \bigg((2^iC_i - y_i)(2^{n-i}C_{n-i}-y_{n-i})\bigg), \textit{ with } y_0=0, \; y_1=1.
\end{equation}}}}
\end{center}

and for large $n$, $\;y_n \sim \Bigg(\frac{10-2\sqrt{10}}{20}\Bigg)\frac{2^{3n-2}}{\sqrt{\pi n^3}}$. 
Where $C_i$ is the $i$th Catalan number.
\end{theorem}
If we look closely to the recurrence relation~(\ref{e:y}), since $d_i= (2^iC_i-y_i)$ 
we obtain the `case 4' in equation~(\ref{eqg:t2}):
\[
y_n =\sum_{i=1}^{n-1} d_i d_{n-i}, \; \textit{ with } \; y_1=1 .
\]
The first ten terms of the sequence $\{y_n\}_{n\geq1}$ are:
\[
1, 1, 6, 29, 162, 978, 6156, 40061, 267338, 819238, \ldots
\]
We have also shown in \cite{P3} that $y_n$ has the following generating function
\begin{center}
\framebox{\parbox[b]{8 cm}{
\[
Y(x) = \frac{2-\sqrt{1-8x} - \sqrt{3-4x-2\sqrt{1-8x}}}{2}.
\]}}
\end{center}
We have also shown in \cite{P3} that the sequence $\{y_n\}_{n\geq 1}$ preserves the parity of Catalan numbers.

\subsubsection{Case 3}

Let $\psi$ and $\phi$ be propositional variables, then
\[\nu(\psi \rightharpoonup  \phi)=1 \quad  : \quad (\nu(\psi) = 1 \; \wedge \; \nu(\phi)=0) .\]

In this case we are interested in formulae obtained from $p_1 \rightharpoonup  \ldots  \rightharpoonup  p_n$ 
by inserting brackets so that the valuation of the first $i$ bracketing give 1 and the rest $(n-i)$ bracketing give 0. 
The table 5 below shows the truth tables, (merged into one), for the two bracketed implications
in $n=3$ variables; where the corresponding case 3 truth values are denoted in green.

\begin{prop}
Let $d_n^{\#3}$ be the number of rows with the value ``true'' in the truth tables 
of all bracketed formulae with $n$ distinct propositions $p_1,\ldots,p_n$ 
connected by the binary connective of m-implication, type 1, such that the valuation 
of the first $i$ bracketing gives 1 and the rest of $(n-i)$ bracketing gives 0. 
Then,
\begin{equation}\label{e:d3}
d_n^{\#3}=\sum_{i=1}^{n-1}d_iy_{n-i}, \textit{ where } \{y_i\} \textit{ is defined in equation~(\ref{e:y}) }.
\end{equation}
Since $d_i=(2^iC_i-y_i)$,
\[
d_n^{\#3}=\sum_{i=1}^{n-1}(2^iC_i-y_i)y_{n-i} .
\]
\end{prop}
\begin{pf}
A row with the value 1, `true', comes from an expression $\psi\rightharpoonup \phi$, where $\nu(\psi)=1$ and $\nu(\phi)=0$. 
If $\psi$ contains $i$ variables, then $\phi$ contains $(n-i)$ variables, and the number of 
choices is given by the summand in the proposition. \qed
\end{pf}

\begin{table}[t]
\caption{$n=3$}
\centering
\begin{tabular}{|c| c | c | c | c|}
\hline \hline 
$p_1$ & $p_2$ & $p_3$ & $p_1\rightharpoonup (p_2 \rightharpoonup p_3)$ & $(p_1 \rightharpoonup p_2)\rightharpoonup p_3$ \\
\hline 1 & 1 & 1 & {\bf \color{green}1} & {\bf \color{red}1} \\
\hline 1 & 1 & 0 & {\bf \color{blue}0} & {\bf 1} \\
\hline 1 & 0 & 1 & {\bf \color{blue}0} & {\bf \color{blue}0} \\
\hline 1 & 0 & 0 & {\bf \color{blue}0} & {\bf \color{green}1} \\
\hline 0 & 1 & 1 & {\bf 1} & {\bf \color{blue}0} \\
\hline 0 & 1 & 0 & {\bf \color{red}1} & {\bf \color{green}1} \\
\hline 0 & 0 & 1 & {\bf \color{red}1} & {\bf \color{blue}0} \\
\hline 0 & 0 & 0 & {\bf \color{red}1} & {\bf \color{green}1} \\
\hline
\end{tabular}
\end{table} 
\medskip 
Let $D_3(x)=\sum_{n\geq 1}d_n^{\#3}x^n$, then $D_3(x) = \sum_{n\geq 1} \sum_{i=1}^{n-1} (2^iC_i-y_i)y_{n-i} x^n$,  
which gives us $D_3(x)= (G(x)-Y(x))Y(x)$, writing this more explicitly gives us the following proposition:

\begin{prop}
The generating function for the sequence $\{ d_n^{\#3} \}_{n\geq 0}$ is given by 
\begin{center}
\framebox{\parbox[b]{15 cm}{
\[
D_3(x)= \frac{-5+3\sqrt{1-8x}+3\sqrt{3-4x-2\sqrt{1-8x}}+4x-\sqrt{1-8x}\sqrt{3-4x-2\sqrt{1-8x}}}{4} .
\]}}
\end{center}

\end{prop}
By using Maple we find the first 25 terms of this sequence:

\begin{eqnarray*}
\{d_n^{\#3} \}_{n\geq 2} &=& 1 , 4 , 19, 108, 646, 4056, 26355, 175628, 1193906, 8246856, \\
&& 57716798, 408391736, 13+2916689516, 20997741104, 152218453443,  \\
&& 1110202813836, 8140864778810, 59981252880360, 443834410644618,  \\
&& 3296876425605992, 24575508928455572, 183773880824034512, \\
&& 1378248141659861486, 10364040821146016568  \ldots
\end{eqnarray*}

\subsubsection{Case 2}

Let $\psi$ and $\phi$ be propositional variables, then
\[\nu(\psi \rightharpoonup  \phi)=1 \quad  : \quad (\nu(\psi) = 0 \; \wedge \; \nu(\phi)=1) .\]

In this case we are interested in formulae obtained from $p_1 \rightharpoonup  \ldots  \rightharpoonup  p_n$ 
by inserting brackets so that the valuation of the first $i$ bracketing gives 0 and the rest $(n-i)$ bracketing give 1. 
The table 5 above shows the truth tables, (merged into one), for the two bracketed implications
in $n=3$ variables; where the corresponding case 3 truth values are denoted in red.

\begin{prop}
Let $d_n^{\#2}$ be the number of rows with the value ``true'' in the truth tables 
of all bracketed formulae with $n$ distinct propositions $p_1,\ldots,p_n$ 
connected by the binary connective of m-implication, type 1, such that the valuation 
of the first $i$ bracketing gives 0 and the rest of $(n-i)$ bracketing gives 1. 
Then,
\begin{equation}\label{e:d2}
d_n^{\#2}=\sum_{i=1}^{n-1}y_id_{n-i}, \textit{ where } \{y_i\}\textit{ is defined in equation~(\ref{e:y})}.
\end{equation}
Since $d_i=(2^iC_i-y_i)$,
\[
d_n^{\#2}=\sum_{i=1}^{n-1}y_i(2^{n-i}C_{n-i}-y_{n-i}) .
\]
\end{prop}
\begin{pf}
A row with the value 1, `true', comes from an expression $\psi\rightharpoonup \phi$, where $\nu(\psi)=0$ and $\nu(\phi)=1$. 
If $\psi$ contains $i$ variables, then $\phi$ contains $(n-i)$ variables, and the number of 
choices is given by the summand in the proposition. \qed
\end{pf}

Let $D_2(x)=\sum_{n\geq 1}d_n^{\#2}x^n$, then $D_2(x) = \sum_{n\geq 1} \sum_{i=1}^{n-1} y_i(2^{n-i}C_{n-i}-y_{n-i})x^n$,  
which gives us $D_2(x)= Y(x)(G(x)-Y(x))$, thus Case 2 coincides with Case 3. Writing this more explicitly gives us the following proposition:

\begin{prop}
The generating function for the sequence $\{ d_n^{\#2} \}_{n\geq 0}$ is given by 
\begin{center}
\framebox{\parbox[b]{15 cm}{
\[
D_2(x)= \frac{-5+3\sqrt{1-8x}+3\sqrt{3-4x-2\sqrt{1-8x}}+4x-\sqrt{1-8x}\sqrt{3-4x-2\sqrt{1-8x}}}{4} .
\]}}
\end{center}

\end{prop}

\subsubsection{Case 1}

Let $\psi$ and $\phi$ be propositional variables, then
\[\nu(\psi \rightharpoonup  \phi)=1 \quad  : \quad ( \nu(\psi) = 0 \; \wedge \; \nu(\phi)=0 ) .\]

In this case we are interested in formulae obtained from $p_1 \rightharpoonup  \ldots  \rightharpoonup  p_n$ 
by inserting brackets so that the valuation of the first $i$ bracketing gives 0 and the rest $(n-i)$ bracketing gives 0. 
The table 5 above shows the truth tables, (merged into one), for the two bracketed implications
in $n=3$ variables; where the corresponding case 3 truth values are denoted in black.

\begin{prop}
Let $d_n^{\#1}$ be the number of rows with the value ``true'' in the truth tables 
of all bracketed formulae with $n$ distinct propositions $p_1,\ldots,p_n$ 
connected by the binary connective of m-implication, type 1, such that the valuation 
of the first $i$ bracketing gives 0 and the rest of $(n-i)$ bracketing gives 0. 
Then,
\begin{equation}\label{e:d1}
d_n^{\#1}=\sum_{i=1}^{n-1}y_iy_{n-i}, \textit{ where } \{y_i\} \textit{ is defined in equation~(\ref{e:y})} .
\end{equation}
\end{prop}
\begin{pf}
A row with the value 1, `true', comes from an expression $\psi\rightharpoonup \phi$, where $\nu(\psi)=0$ and $\nu(\phi)=0$. 
If $\psi$ contains $i$ variables, then $\phi$ contains $(n-i)$ variables, and the number of 
choices is given by the summand in the proposition. \qed
\end{pf}

Let $D_1(x)=\sum_{n\geq 1}d_n^{\#1}x^n$, then $D_1(x) = \sum_{n\geq 1} \sum_{i=1}^{n-1}y_iy_{n-i} x^n$,  
which gives us $D_2(x)= Y(x)^2$. Writing this more explicitly gives us the following proposition:

\begin{prop}
The generating function for the sequence $\{ d_n^{\#1} \}_{n\geq 0}$ is given by 

\begin{center}
\framebox{\parbox[b]{15 cm}{
\[
D_1(x)= \frac{4-3\sqrt{1-8x}-2\sqrt{3-4x-2\sqrt{1-8x}}-6x+\sqrt{1-8x}\sqrt{3-4x-2\sqrt{1-8x}}}{2} .
\]}}
\end{center}
\end{prop}

By using Maple we find the first 25 terms of this sequence:

\begin{eqnarray*}
\{d_n^{\#1} \}_{n\geq 2} &=& 1, 2, 13, 70, 418, 2628, 17053, 113566, 771638, 5327804, 37274482, 263669500, \\
&& 1882630692, 13550468360, 98212733277, 716195167502, 5250931034798,   \\
&& 8683418448780, 286206574421222, 2125766544922612, 15844332066531484,  \\
&& 3118472460044221368, 888436633672089842, 6680306733514013388, \ldots
\end{eqnarray*}

\subsection{Type 2}
\begin{defn}
Let $\psi$, and $\phi$ be propositional variables then 
\[
\psi \leftharpoonup \phi \equiv \neg\psi \to \phi .
\]
For any valuation $\nu$,
\[\nu(\psi\rightharpoonup\phi)=\cases{0 & if $\nu(\psi)=0$ and $\nu(\phi)=0$,\cr
1 & otherwise.\cr}\] 
\end{defn}

\medskip

\medskip

Let  $k_n$, $h_n$ be the number of rows with the value ``true'', and ``false'' in the truth tables 
of all bracketed formulae with $n$ distinct propositions $p_1,\ldots,p_n$ 
connected by the binary connective of m-implication, in the type (ii), respectively. Then using theorem \ref{T:1}, 
\[
g_n = \sum_{i=1}^{n-1} (h_i+k_i)(h_{n-1}+k_{n-i})
\]
Thus if we expand the right hand side 
\begin{equation} \label{eq2.1}
g_n=\underbrace{\sum_{i=1}^{n-1}h_i h_{n-i}}_{\textit{case 4}} + \underbrace{\sum_{i=1}^{n-1}h_i k_{n-i}}_{\textit{case 3}}  + \underbrace{\sum_{i=1}^{n-1}k_i h_{n-i}}_{\textit{case 2}}   + \underbrace{\sum_{i=1}^{n-1}k_i k_{n-i}}_{\textit{case 1}}  .
\end{equation}
Let
\[
k_n^{\#1}= \sum_{i=1}^{n-1}k_ik_{n-i}, \quad k_n^{\#2}= \sum_{i=1}^{n-1}k_ih_{n-i}, \quad k_n^{\#3}= \sum_{i=1}^{n-1}h_ik_{n-i},\quad h_n= \sum_{i=1}^{n-1}h_ih_{n-i} .   
\]
Where
\[
0=k_1^{\#1}=k_1^{\#2}=k_1^{\#3}, \quad \textit{ and }\quad h_1=1 .
\]

\subsubsection{Case 4}
This case has already been discussed in \cite{P3}, and we had the following results:
\begin{prop}
Let $h_n$ be the number of rows with the value ``false'' in the truth tables 
of all bracketed formulae with $n$ distinct propositions $p_1,\ldots,p_n$ 
connected by the binary connective of m-implication, type 2, such that the valuation 
of the first $i$ bracketing gives 0 and the rest of $(n-i)$ bracketing gives 0. 
Then,
\begin{equation}\label{e:h}
h_n=\sum_{i=1}^{n-1} h_ih_{n-i}, \textit{ where } h_0=0, h_1=1.
\end{equation}
\end{prop}
The recurrence relation (\ref{e:h}) is very well known; it is the recurrence relation for Catalan numbers. The Catalan numbers has the following generating function, and the explicit formula:
\[
H(x)=\frac{1-\sqrt{1-4x}}{2}, \quad h_n=\frac{1}{n}{{2n-2}\choose{n-1}} .
\]

\begin{cor}
Suppose  we have all possible well-formed formulae obtained from
$p_1 \leftharpoonup p_2 \leftharpoonup \ldots \leftharpoonup p_n$ by inserting brackets, 
where $p_1,\ldots,p_n$ are distinct propositions. 
Then each formula defines the same truth table, (provided that we do not distinguish the 1s by their case type).
\end{cor}

 Here are the truth tables, (merged into one), for the bracketed m-implications, in $n=3$ variables. Where the corresponding 
false entries are denoted in blue:

\begin{table}[h]
\caption{$n=3$}
\centering
\begin{tabular}{|c| c | c | c | c|}
\hline \hline 
\hline $p_1$ & $p_2$ & $p_3$ & $p_1\leftharpoonup (p_2 \leftharpoonup p_3)$ & $(p_1 \leftharpoonup p_2)\leftharpoonup p_3$ \\
\hline 1 & 1 & 1 &  {\bf \color{red}1} & {\bf 1} \\
\hline 1 & 1 & 0 &  {\bf \color{red}1}  &   {\bf \color{green}1} \\
\hline 1 & 0 & 1 &  {\bf \color{red}1}  & {\bf 1} \\
\hline 1 & 0 & 0 &  {\bf \color{green}1}  &  {\bf \color{green}1} \\
\hline 0 & 1 & 1 & {\bf 1} & {\bf 1} \\
\hline 0 & 1 & 0 & {\bf 1} &  {\bf \color{green}1} \\
\hline 0 & 0 & 1 & {\bf 1} &  {\bf \color{red}1}  \\
\hline 0 & 0 & 0 & {\bf \color{blue} 0} & {\bf \color{blue} 0} \\
\hline
\end{tabular}
\end{table} 
\medskip

\subsubsection{Case 3}

Let $\psi$ and $\phi$ be propositional variables, then
\[\nu(\psi \leftharpoonup  \phi)=1 \quad  : \quad (\nu(\psi) = 0 \; \wedge \; \nu(\phi)=1) .\]

In this case we are interested in formulae obtained from $p_1 \leftharpoonup p_2 \leftharpoonup \ldots \leftharpoonup p_n$ 
by inserting brackets so that the valuation of the first $i$ bracketing gives 0 and the rest $(n-i)$ bracketing gives 0. 
The table 6 above shows the truth tables, (merged into one), for the two bracketed implications
in $n=3$ variables. The corresponding case 3 truth values are denoted in red.

\begin{prop}
Let $k_n^{\#3}$ be the number of rows with the value ``true'' in the truth tables 
of all bracketed formulae with $n$ distinct propositions $p_1,\ldots,p_n$ 
connected by the binary connective of m-implication, type 2, such that the valuation 
of the first $i$ bracketing gives 0 and the rest of $(n-i)$ bracketing gives 1. 
Then,
\begin{equation}\label{e:t2k3}
k_n^{\#3}=\sum_{i=1}^{n-1}h_ik_{n-i}, \textit{ where } \{h_i\} \textit{ is defined in equation~(\ref{e:h})}.
\end{equation}
\end{prop}
\begin{pf}
A row with the value 1, `true', comes from an expression $\psi\rightharpoonup \phi$, where $\nu(\psi)=0$ and $\nu(\phi)=1$. 
If $\psi$ contains $i$ variables, then $\phi$ contains $(n-i)$ variables, and the number of 
choices is given by the summand in the proposition. \qed
\end{pf}
Since $k_i=(2^iC_i-h_i)$, from equation~(\ref{e:t2k3}) we get,
\[
k_n^{\#3}= \sum_{i=1}^{n-1} h_i(2^{n-i}C_{n-i}-h_{n-i}). 
\]
Let $K_3(x)=\sum_{i\geq 1} k_n^{\#3} x^n$, then we get 
$K_3(x)= H(x)(G(x)-H(x))$, writing the right hand side more explicitly gives the following proposition:

\begin{prop}
The generating function for the sequence $\{ k_n^{\#3} \}_{n\geq 0}$ is given by 
\begin{center}
\framebox{\parbox[b]{11 cm}{
\[
K_3(x)= \frac{-1-\sqrt{1-8x}+\sqrt{1-8x}\sqrt{1-4x}+\sqrt{1-4x}+4x}{4} .
\]}}
\end{center}

\end{prop}
By using Maple we find the first 25 terms of this sequence:

\begin{eqnarray*}
\{k_n^{\#3} \}_{n\geq 2} &=& 1, 4, 19, 100, 566, 3384, 21107, 136084, 900674, 6087496, \\
&&41850366, 291766952, 2057964492, 14659421040, 105305580483, 761981900724, \\
&& 5548736343434, 0632122219688, 299017702596554, 2210275626304248, \\
&& 16403005547059508, 122169144755555088, 912887876722311406, 684174390763667239, \ldots
\end{eqnarray*}

\subsubsection{Case 2}

Let $\psi$ and $\phi$ be propositional variables, then
\[\nu(\psi \leftharpoonup  \phi)=1 \quad  : \quad (\nu(\psi) = 1 \; \wedge \; \nu(\phi)=0) .\]

In this case we are interested in formulae obtained from $p_1 \leftharpoonup p_2 \leftharpoonup \ldots \leftharpoonup p_n$ 
by inserting brackets so that the valuation of the first $i$ bracketing gives 1 and the rest $(n-i)$ bracketing gives 0. 
The table 6 above shows the truth tables, (merged into one), for the two bracketed implications
in $n=3$ variables. The corresponding case 3 truth values are denoted in green.

\begin{prop}
Let $k_n^{\#2}$ be the number of rows with the value ``true'' in the truth tables 
of all bracketed formulae with $n$ distinct propositions $p_1,\ldots,p_n$ 
connected by the binary connective of m-implication, type 2, such that the valuation 
of the first $i$ bracketing gives 1 and the rest of $(n-i)$ bracketing gives 0. 
Then,
\begin{equation}\label{e:t2k2}
k_n^{\#2}=\sum_{i=1}^{n-1}k_ih_{n-i}, \textit{ where } \{h_i\} \textit{ is defined in equation ~(\ref{e:h})}.
\end{equation}
\end{prop}
\begin{pf}
A row with the value 1, `true', comes from an expression $\psi\rightharpoonup \phi$, where $\nu(\psi)=1$ and $\nu(\phi)=0$. 
If $\psi$ contains $i$ variables, then $\phi$ contains $(n-i)$ variables, and the number of 
choices is given by the summand in the proposition. \qed
\end{pf}
Since $k_i=(2^iC_i-h_i)$, from equation~(\ref{e:t2k2}) we get,
\[
k_n^{\#2}= \sum_{i=1}^{n-1} (2^{i}C_{i}-h_{i})h_{n-i}. 
\]
Let $K_2(x)=\sum_{i\geq 1} k_n^{\#2} x^n$, then we get 
$K_2(x)= (G(x)-H(x))H(x)$, writing the right hand side more explicitly gives the following proposition:

\begin{prop}
The generating function for the sequence $\{ k_n^{\#2} \}_{n\geq 0}$ is given by 
\begin{center}
\framebox{\parbox[b]{11 cm}{
\[
K_2(x)= \frac{-1-\sqrt{1-8x}+\sqrt{1-8x}\sqrt{1-4x}+\sqrt{1-4x}+4x}{4} .
\]
}}
\end{center}
\end{prop}

\subsubsection{Case 1}

Let $\psi$ and $\phi$ be propositional variables, then
\[\nu(\psi \leftharpoonup  \phi)=1 \quad  : \quad (\nu(\psi) = 1 \; \wedge \; \nu(\phi)=1) .\]

In this case we are interested in formulae obtained from $p_1 \leftharpoonup p_2 \leftharpoonup \ldots \leftharpoonup p_n$ 
by inserting brackets so that the valuation of the first $i$ bracketing gives 1 and the rest $(n-i)$ bracketing gives 1. 
The table 6 above shows the truth tables, (merged into one), for the two bracketed implications
in $n=3$ variables. The corresponding case 3 truth values are denoted in black.
\begin{prop}
Let $k_n^{\#1}$ be the number of rows with the value ``true'' in the truth tables 
of all bracketed formulae with $n$ distinct propositions $p_1,\ldots,p_n$ 
connected by the binary connective of m-implication, type 2, such that the valuation 
of the first $i$ bracketing gives 1 and the rest of $(n-i)$ bracketing gives 1. 
Then,
\begin{equation}\label{e:t2k1}
k_n^{\#1}=\sum_{i=1}^{n-1}k_ik_{n-i}.
\end{equation}
\end{prop}
\begin{pf}
A row with the value 1, `true', comes from an expression $\psi\rightharpoonup \phi$, where $\nu(\psi)=1$ and $\nu(\phi)=1$. 
If $\psi$ contains $i$ variables, then $\phi$ contains $(n-i)$ variables, and the number of 
choices is given by the summand in the proposition. \qed
\end{pf}
Since $k_i=(2^iC_i-h_i)$, from equation~(\ref{e:t2k3}) we get,
\[
k_n^{\#2}= \sum_{i=1}^{n-1} (2^iC_i-h_i)(2^{n-i}C_{n-i}-h_{n-i}), \quad \textit{ where $h_i$ is the sequence in~(\ref{e:h})}.
\]
Let $K_1(x)=\sum_{i\geq 1} k_n^{\#1} x^n$, then we get 
$K_2(x)= (G(x)-H(x))^2$, writing the right hand side more explicitly gives the following proposition:

\begin{prop}
The generating function for the sequence $\{ k_n^{\#1} \}_{n\geq 0}$ is given by 

\begin{center}
\framebox{\parbox[b]{7 cm}{
\[
K_1(x)= \frac{1-6x-\sqrt{1-4x}\sqrt{1-8x}}{2} .
\]
}}
\end{center}
\end{prop}
By using Maple we find the first 25 terms of this sequence:

\begin{eqnarray*}
\{k_n^{\#1} \}_{n\geq 2} &=& 1, 6, 37, 234, 514, 9996, 67181, 458562, 3172478, 22206420, \\
&& 157027938, 1120292388, 8055001716, 58314533400, 424740506109, 3110401363122,  \\
&& 22888001498102, 169155516667524, 1255072594261142, 9345400450314924, \\
&&  69812926066668044, 523072984217339304, 3929809142578361938, 29598511892723647860, \ldots
\end{eqnarray*}

\subsection{Type 3}

\begin{defn}
Let $\psi$, and $\phi$ be propositional variables then 
\[
\psi \rightleftharpoons \phi \equiv \neg\psi \to \neg\phi .
\]
For any valuation $\nu$,
\[\nu(\psi \rightleftharpoons\phi)=\cases{0 & if $\nu(\psi)=0$ and $\nu(\phi)=1$,\cr
1 & otherwise.\cr}\] 
\end{defn}

\medskip

\medskip

Let  $b_n$, $s_n$ be the number of rows with the value ``true'', and ``false'' in the truth tables 
of all bracketed formulae with $n$ distinct propositions $p_1,\ldots,p_n$ 
connected by the binary connective of m-implication, in the type (iii), respectively. Then using theorem \ref{T:1}, 
\[
g_n = \sum_{i=1}^{n-1} (b_i+s_i)(b_{n-i}+s_{n-i})
\]
Thus if we expand the right hand side 
\begin{equation} \label{eq2.2}
g_n=\underbrace{\sum_{i=1}^{n-1}b_i b_{n-i}}_{\textit{case 1}} + \underbrace{\sum_{i=1}^{n-1} b_i s_{n-i}}_{\textit{case 2}}  + \underbrace{\sum_{i=1}^{n-1}s_i b_{n-i}}_{\textit{case 4}}   + \underbrace{\sum_{i=1}^{n-1}s_i s_{n-i}}_{\textit{case 3}}  .
\end{equation}
Let
\[
b_n^{\#1}= \sum_{i=1}^{n-1}b_ib_{n-i}, \quad b_n^{\#2}= \sum_{i=1}^{n-1}b_is_{n-i}, \quad b_n^{\#3}= \sum_{i=1}^{n-1}s_is_{n-i},\quad s_n= \sum_{i=1}^{n-1}s_ib_{n-i} .   
\]
Where
\[
0=b_1^{\#1}=b_1^{\#2}=b_1^{\#3}, \quad \textit{ and }\quad s_1=1 .
\]

\begin{prop}
Let $s_n^{\#3}$ be the number of rows with the value ``false'' in the truth tables 
of all bracketed formulae with $n$ distinct propositions $p_1,\ldots,p_n$ 
connected by the binary connective of m-implication, type 3, such that the valuation 
of the first $i$ bracketing gives 0 and the rest of $(n-i)$ bracketing gives 1. 
Then,
\begin{equation}\label{e:s}
s_n = \sum_{i=1}^{n-1} s_i b_{n-i}, \textit{ where } s_0=0, s_1=1.
\end{equation}
Since $b_i=2^iC_i-s_i$, 
\[
s_n = \sum_{i=1}^{n-1} s_i (2^{n-i}C_{n-i}-s_{n-i}),
\]
\end{prop}
\begin{pf}
A row with the value 1, `true', comes from an expression $\psi\rightleftharpoons \phi$, where $\nu(\psi)=0$ and $\nu(\phi)=1$. 
If $\psi$ contains $i$ variables, then $\phi$ contains $(n-i)$ variables, and the number of 
choices is given by the summand in the proposition. \qed
\end{pf}

Since the recurrence relation~(\ref{e:s}), is equivalent to the recurrence relation for the sequence 
$\{f_n\}_{n\geq 1}$, the four generating function in section 3, (starting from page 5), and the corresponding asymptotics, (starting from page 11), 
will appear again. Thus the reader is suggested to work these cases on his-own.

\section{Asymptotic Estimate 2}

\subsection{Asymptotic for Type 1 and 3}

In \cite{P3} we have shown that the following asymptotic results are true for the sequence $\{y_n\}_{n \geq 1}$ :

\begin{theorem}
Let $y_n$ be number of rows with the value false in the truth tables of all the
bracketed m-implications, case(i), with $n$ distinct variables. Then
\begin{center}
\framebox{\parbox[b]{6 cm}{
\[y_n \sim \left(\frac{10-2\sqrt{10}}{10}\right)\frac{2^{3n-2}}{\sqrt{\pi n^3}}.\]
}}
\end{center}

\end{theorem}

\begin{cor}
Let $g_n$ be the total number of rows in all truth tables for bracketed
m-implications, case(i), with $n$ distinct variables, and $y_n$ the number of rows with the
value ``false''. Then $\lim_{n\to\infty}y_n/g_n= \frac{10-2\sqrt{10}}{10}$.
\end{cor}
 
\begin{theorem}
Let $d_n^{\#3}$ be number of rows with the value true in the truth tables of all the
bracketed implications with $n$ variables, type 1, case 3. Then

\begin{center}
\framebox{\parbox[b]{6 cm}{
\[ d_n^{\#3} \sim \left(\frac{11\sqrt{10}-30}{20}\right)
\frac{2^{3n-2}}{\sqrt{\pi n^3}}.\]
}}
\end{center}

\end{theorem}
\begin{pf}
Recall that 

{\small
\[
D_3(x)= \frac{-5+3\sqrt{1-8x}+3\sqrt{3-4x-2\sqrt{1-8x}}+4x-\sqrt{1-8x}\sqrt{3-4x-2\sqrt{1-8x}}}{4} .
\]
}
By using the asymptotic techniques which we have discussed in \cite{P1} and \cite{P3}, we found that
$r=\frac{1}{8}$, $f(x)= \sqrt{1-8x}$ and $D_3(\frac{1}{8})= \frac{-9+3\sqrt{10}}{8}\not= 0$. Let
$A(x)=D_3(x)-D_3(\frac{1}{8})$

\[\lim_{x\to1/8} \frac{A(x)}{f(x)} 
= \lim_{x\to 1/8} \frac{-1+6\sqrt{1-8x}+6\sqrt{2-4x-2\sqrt{1-8x}}-2\sqrt{1-8x}\sqrt{2-4x-2\sqrt{1-8x}}+8x-3\sqrt{10}}{8\sqrt{1-8x}} .\]

Let $v=\sqrt{1-8x}$. Then
\begin{eqnarray*}
L &=& \lim_{v\to 0} \frac{6v+3\sqrt{2}\sqrt{5-4v+v^2}-\sqrt{2}v\sqrt{5-4v+v^2}-v^2-3\sqrt{10}}{8v} \\
&=& \lim_{v\to 0} -\frac{-6\sqrt{5-4v+v^2}+11\sqrt{2}-9\sqrt{2}v^2+2v\sqrt{5-4v-v^2}}{8\sqrt{5-4v-v^2}} \\
&=& - \frac{\sqrt{5}(-6\sqrt{5}+11\sqrt{2})}{40}\\
&=& - \frac{11\sqrt{10}-30}{40},
\end{eqnarray*}
where we have used l'H\^opital's Rule in the penultimate line.

Finally,
\[t_n^{\#3} \sim - \frac{11\sqrt{10}-30}{40}{n-\frac{3}{2}\choose n}
\left(\frac{1}{8}\right)^{-n} \sim \left(\frac{11\sqrt{10}-30}{20}\right)
\frac{2^{3n-2}}{\sqrt{\pi n^3}},\]
and the proof is finished. \qed
\end{pf}

\begin{theorem}
Let $d_n^{\#2}$ be number of rows with the value true in the truth tables of all the
bracketed implications with $n$ variables, type 1, case 2. Then
\begin{center}
\framebox{\parbox[b]{6 cm}{
\[ d_n^{\#2} \sim \left(\frac{11\sqrt{10}-30}{20}\right)
\frac{2^{3n-2}}{\sqrt{\pi n^3}}.\]
}}
\end{center}

\end{theorem}

\begin{theorem}
Let $d_n^{\#1}$ be number of rows with the value true in the truth tables of all the
bracketed implications with $n$ variables, type 1, case 1. Then

\begin{center}
\framebox{\parbox[b]{6 cm}{
\[ d_n^{\#1} \sim \left(\frac{20-9\sqrt{10}}{10}\right)
\frac{2^{3n-2}}{\sqrt{\pi n^3}}.\]

}}
\end{center}

\end{theorem}

\begin{pf}
Recall that 

{\small
\[
D_1(x)= \frac{4-3\sqrt{1-8x}-2\sqrt{3-4x-2\sqrt{1-8x}}-6x+\sqrt{1-8x}\sqrt{3-4x-2\sqrt{1-8x}}}{2} .
\]
}
By using the asymptotic techniques which we have discussed in \cite{P1} and \cite{P3}, we found that
$r=\frac{1}{8}$, $f(x)= \sqrt{1-8x}$ and $D_1(\frac{1}{8})= \frac{13-4\sqrt{10}}{8}\not= 0$. Let
$A(x)=D_1(x)-D_1(\frac{1}{8})$

\[\lim_{x\to1/8} \frac{A(x)}{f(x)} 
= \lim_{x\to 1/8} \frac{3-12\sqrt{1-8x}-8\sqrt{2-4x-2\sqrt{1-8x}}-24x +4\sqrt{1-8x}\sqrt{2-4x-2\sqrt{1-8x}} +4\sqrt{10}}{8\sqrt{1-8x}}  .\]

Let $v=\sqrt{1-8x}$. Then
\begin{eqnarray*}
L &=& \lim_{v\to 0} \frac{-12v-4\sqrt{2}\sqrt{v^2-4v+5}+3v^2+2\sqrt{2}v\sqrt{v^2-4v+5}+4\sqrt{10}}{8v} \\
&=& \lim_{v\to 0} -\frac{1}{4}\frac{-6\sqrt{v^2-4v+5}-8\sqrt{2}v+9\sqrt{2}+3v\sqrt{v^2-4v+5}+2\sqrt{2}x^2}{\sqrt{v^2-4v+5}} \\
&=& - \frac{\sqrt{5}(-12\sqrt{5}+18\sqrt{2})}{40}\\
&=& - \frac{30-9\sqrt{10}}{20},
\end{eqnarray*}
where we have used l'H\^opital's Rule in the penultimate line.

Finally,
\[d_n^{\#1} \sim - \frac{30-9\sqrt{10}}{20}{n-\frac{3}{2}\choose n}
\left(\frac{1}{8}\right)^{-n} \sim \left(\frac{30-9\sqrt{10}}{10}\right)
\frac{2^{3n-2}}{\sqrt{\pi n^3}},\]
and the proof is finished. \qed
\end{pf}

\begin{cor}
Let $g_n$ be the total number of rows in all truth tables for bracketed
implications with $n$ variables, then 

\[
\lim_{n\to\infty}\frac{d_n^{\#1}}{g_n}=\frac{30-9\sqrt{10}}{10}, \;\;\;  \lim_{n\to\infty}\frac{d_n^{\#3}}{g_n}=\lim_{n\to\infty}\frac{d_n^{\# 2}}{g_n}=\frac{11\sqrt{10}-30}{20}, \;\;\; \lim_{n\to\infty}\frac{y_n}{g_n}=\frac{10-2\sqrt{10}}{10}. 
\]
\end{cor}

\begin{cor}
\[
\lim_{n\to\infty}\frac{d_n^{\#1}}{d_n^{\#3}}=\frac{12\sqrt{10}-18}{31}, \;\;\; \lim_{n\to\infty}\frac{d_n^{\#3}}{d_n^{\#1}}=\frac{2\sqrt{10}+3}{6}, \;\;\; \lim_{n\to\infty}\frac{d_n^{\#1}}{y}=\frac{4-\sqrt{10}}{2}, 
\]
\[
\lim_{n\to\infty}\frac{y}{d_n^{\#1}}=\frac{4+\sqrt{10}}{3}, \;\;\; \lim_{n\to\infty}\frac{y}{d_n^{\#3}}=\frac{16+10\sqrt{10}}{31}, \;\;\; \lim_{n\to\infty}\frac{d_n^{\#3}}{y}=\frac{5\sqrt{10}-8}{12}. 
\]

\end{cor}

\subsection{Asymptotic for Type 2}

In \cite{P3} we have shown that the following asymptotic results are true for the sequence $\{h_n\}_{n\geq 1}$:

\begin{theorem}
Let $h_n$ be number of rows with the value false in the truth tables of all the
bracketed m-implications, case(i), with $n$ distinct variables. Then
\begin{center}
\framebox{\parbox[b]{4 cm}{
\[h_n \sim \frac{2^{2n}}{\sqrt{\pi n^3}}.\]
}}
\end{center}

\end{theorem}

\begin{cor}
Let $g_n$ be the total number of rows in all truth tables for bracketed
m-implications, case(i), with $n$ distinct variables, and $y_n$ the number of rows with the
value ``false', in type 2'. Then $\lim_{n\to\infty}h_n/g_n= 0$.
\end{cor}
\begin{pf}
\[
\lim_{n \to \infty} \frac{h_n}{g_n} = \lim_{n\to \infty } \frac{1}{2^{n-2}} = 0 . 
\]
\qed
\end{pf}

\begin{theorem}
Let $k_n^{\#1}$ be number of rows with the value true in the truth tables of all the
bracketed implications with $n$ variables, type 1, case 1. Then

\begin{center}
\framebox{\parbox[b]{6 cm}{
\[ k_n^{\#3} \sim \left(\frac{\sqrt{2}}{2}\right)
\frac{2^{3n-2}}{\sqrt{\pi n^3}}.\]
}}
\end{center}

\end{theorem}
\begin{pf}
Recall that 

\[
K_1(x)= \frac{1-6x-\sqrt{1-4x}\sqrt{1-8x}}{2} .
\]

By using the asymptotic techniques which we have discussed in \cite{P1} and \cite{P3}, we found that
$r=\frac{1}{8}$, $f(x)= \sqrt{1-8x}$ and $K_1(\frac{1}{8})= \frac{1}{8}\not= 0$. Let
$A(x)=K_1(x)-K_1(\frac{1}{8})$

\[\lim_{x\to1/8} \frac{A(x)}{f(x)} 
= \lim_{x\to 1/8} \frac{3-24x-4\sqrt{1-8x}\sqrt{1-4x}}{8\sqrt{1-8x}}  .\]

Let $v=\sqrt{1-8x}$. Then
\begin{eqnarray*}
L &=& \lim_{v\to 0} \frac{3v^2-2\sqrt{2}v\sqrt{1+v^2}}{8v} \\
&=& \lim_{v\to 0} -\frac{2(-3v\sqrt{1+v^2}+\sqrt{2}+2\sqrt{2}v^2)}{1+\sqrt{v^2}} \\
&=& - \frac{\sqrt{2}}{4},
\end{eqnarray*}
where we have used l'H\^opital's Rule in the penultimate line.

Finally,
\[k_n^{\#1} \sim - \frac{\sqrt{2}}{4}{n-\frac{3}{2}\choose n}
\left(\frac{1}{8}\right)^{-n} \sim \left(\frac{\sqrt{2}}{2}\right)
\frac{2^{3n-2}}{\sqrt{\pi n^3}},\]
and the proof is finished. \qed
\end{pf}

\begin{cor}
Let $g_n$ be the total number of rows in all truth tables for bracketed
m-implications, case(i), with $n$ distinct variables, and $k_n^{\#1}$ be 
the number of rows with the value `true', in type 2. Then $\lim_{n\to\infty}k_n^{\#1}/g_n= \frac{\sqrt{2}}{2}$.
\end{cor}

\begin{theorem}
Let $k_n^{\#3,2}$ be number of rows with the value true in the truth tables of all the
bracketed implications with $n$ variables, type 1, case 1. Then

\begin{center}
\framebox{\parbox[b]{6 cm}{
\[ k_n^{\#3,2} \sim \left(\frac{2-\sqrt{2}}{4}\right)
\frac{2^{3n-2}}{\sqrt{\pi n^3}}.\]
}}
\end{center}

\end{theorem}

\begin{pf}
Since the total probability has to add up to 1, and 
$K_2(x)$ and $K_3(x)$ both have the same generating function:
\[
\frac{1-\frac{\sqrt{2}}{2}}{2}= \frac{2-\sqrt{2}}{4} .
\]
\qed
\end{pf}

\begin{cor}
Let $g_n$ be the total number of rows in all truth tables for bracketed
implications with $n$ variables, then 
\[
\lim_{n\to\infty}\frac{k_n^{\#2}}{g_n}= \lim_{n\to\infty}\frac{k_n^{\#3}}{g_n}=\frac{2-\sqrt{2}}{4}, \;\;\;  \lim_{n\to\infty}\frac{k_n^{\#1}}{g_n}=\frac{\sqrt{2}}{2}, \;\;\; \lim_{n\to\infty}\frac{h_n}{g_n}= 0. 
\]
\end{cor}

\begin{cor}
\[
\lim_{n\to\infty}\frac{k_n^{\#2}}{k_n^{\#1}}=\frac{\sqrt{2}-1}{2}, 
\;\;\; \lim_{n\to\infty}\frac{k_n^{\#1}}{k_n^{\#2}}=2\sqrt{2}+2.
\]
\end{cor}

\section{Parity}
Recall the recurrence relation of $g_n$
\[
g_n =\sum_{i=1}^{n-1}g_ig_{n-i}
\]
Let $a_i$ and $b_i$ be the corresponding number of false, and truth entries in the considered truth table then:

\begin{equation}
 g_n=\sum_{i=1}^{n-1} (a_i+b_i)(a_{n-i}+b_{n-i}) =\sum_{i=1}^{n-1}a_i a_{n-i} + \sum_{i=1}^{n-1}a_i b_{n-i} + \sum_{i=1}^{n-1}b_i a_{n-i} + \sum_{i=1}^{n-1}b_i b_{n-i}
\end{equation}
Let
\[
\quad z_n^{\#1}=\sum_{i=1}^{n-1}a_i a_{n-i}, \quad  z_n^{\#2}=\sum_{i=1}^{n-1}a_i b_{n-i},\quad   z_n^{\#3}=\sum_{i=1}^{n-1}b_i a_{n-i},\quad   z_n^{\#4}=\sum_{i=1}^{n-1}b_i b_{n-i} .
\]

\begin{theorem}\label{f:c}For $i=1,2,3,4$ each sequence $z_n^{\#i}$ preserves the parity of Catalan numbers.
\end{theorem}

\begin{pf} 
Since $a_i=(2^iC_i-b_i)$ then
\[
\quad 
z_n^{\#1}=\sum_{i=1}^{n-1}(2^iC_i-b_i)(2^{n-i}C_{n-i}-b_{n-i}), \quad  z_n^{\#2}=\sum_{i=1}^{n-1} (2^iC_i-b_i)b_{n-i},\quad   
\]
\[
z_n^{\#3}=\sum_{i=1}^{n-1}b_i (2^{n-i}C_{n-i}-b_{n-i}),\quad   z_n^{\#4}=\sum_{i=1}^{n-1}b_i b_{n-i} .
\]
If an additive partition of $z_n^{\#1}$, is odd, then it comes as a pair; i.e.  
\[
(2^iC_i-b_i)(2^{n-i}C_{n-i}-b_{n-i}) \in \mathbb{O} \Longleftrightarrow b_i, b_{n-i}\in\mathbb{O}.
\]
Hence, $\bigg((2^iC_i-b_i)(2^{n-i}C_{n-i}-b_{n-i})+(2^{n-i}C_{n-i}-b_{n-i})(2^iC_i-b_i)\bigg) \in \mathbb{E}.$\\\\
Thus, $z_n^{\#1}$ can be expressed as a piecewise function depending on the parity of $n$:
\[z_n^{\#1}=\cases{2\sum_{i=1}^{\frac{n-1}{2}} ((2^iC_i-b_i)(2^{n-i}C_{n-i}-b_{n-i}) ) \;\;& if $n\in \mathbb{O}$,\cr\cr
\bigg(2\sum_{i=1}^{\frac{n-2}{2}} ((2^iC_i-b_i)(2^{n-i}C_{n-i}-b_{n-i}) ) \bigg)+ (2^{\frac{n}{2}}C_{\frac{n}{2}}-b_{\frac{n}{2}})^2\;\; & if $n\in \mathbb{E}$.\cr}\]
Finally,  
\[
z_n^{\#1}\in \mathbb{O} \Longleftrightarrow (2^{\frac{n}{2}}C_{\frac{n}{2}}-b_{\frac{n}{2}})^2 \in \mathbb{O} \Longleftrightarrow b_\frac{n}{2} \in \mathbb{O}  \Longleftrightarrow  n=2^i, \;\; \forall i\in\mathbb{N}.
\]
\medskip

Similarly, if an additive partition of $z_n^{\#2}$, is odd, then it comes as a pair; i.e.  
\[
(2^iC_i-b_i)b_{n-i} \in \mathbb{O} \Longleftrightarrow b_i, b_{n-i}\in \mathbb{O} \Longleftrightarrow (2^{n-i}C_{n-i}-b_{n-i})b_i\in \mathbb{O} .
\]
Hence, $\bigg((2^iC_i-b_i)f_{n-i}+(2^{n-i}C_{n-i}-b_{n-i})b_i \bigg) \in \mathbb{E}.$\\\
Thus, $z_n^{\#2}n$ can be expressed as a piecewise function depending on the parity of $n$:
{\small
\[z_n^{\#2}=\cases{\sum_{i=1}^{\frac{n-1}{2}} ((2^iC_i-b_i)b_{n-i}+(2^{n-i}C_{n-i}-b_{n-i})b_i ) & if $n\in \mathbb{O}$,\cr\cr
\bigg(\sum_{i=1}^{\frac{n-2}{2}} ((2^iC_i-b_i)b_{n-i}+(2^{n-i}C_{n-i}-b_{n-i})b_i )\bigg)
+ (2^{\frac{n}{2}}C_{\frac{n}{2}}-b_{\frac{n}{2}})b_{\frac{n}{2}}  & if $n\in \mathbb{E}$.\cr}\]
}
Finally,  
\[
z_n^{\#2}\in \mathbb{O} \Longleftrightarrow  (2^{\frac{n}{2}}C_{\frac{n}{2}}-b_{\frac{n}{2}})b_{\frac{n}{2}} \in \mathbb{O} \Longleftrightarrow b_\frac{n}{2} \in \mathbb{O}  \Longleftrightarrow  n=2^i, \;\; \forall i\in\mathbb{N}.
\]
\\\\
The sequences $z_n^{\#3}$ preserves the parity of Catalan numbers, observe that 
\[
z_n^{\#3}=\sum_{i=1}^{n-1}b_i(2^{n-1}C_{n-i}-b_{n-i}) = \sum_{i=1}^{n-1} (2^iC_i-b_i)b_{n-i} = z_n^{\#2}.
\]
The sequence $z_n^{\#4}$ preserves the parity of the Catalan numbers,  since
\[z_n^{\#4}=\cases{2(b_1b_{n-1}+b_2b_{n-2}+\ldots+b_{\frac{n-1}{2}}b_{\frac{n+1}{2}}) \;\;& if $n\in \mathbb{O}$,\cr\cr
2(b_1b_{n-1}+b_2b_{n-2}+\ldots+b_{\frac{n-2}{2}}b_{\frac{n+2}{2}})+b_{\frac{n}{2}}^2 \;\; & if $n\in \mathbb{E}$.\cr}\]
For $n\geq 2$ 
\[
z_n^{\#4}\in \mathbb{O} \Longleftrightarrow b_{\frac{n}{2}}^2 \in \mathbb{O}  \Longleftrightarrow 
b_{\frac{n}{2}}\in \mathbb{O} \Longleftrightarrow n=2^i \;\;\forall i\in \mathbb{N}. 
\]
\qed
\end{pf}

\pagebreak

\begin{table}[b]
\begin{verse}
Ne yaz\i k, pi\c{s}mi\c{s} ekmek \c{c}i\~glerin elinde;\\
Ne yaz\i k, \c{c}e\c{s}meler cimrilerin elinde.\\
O can\i m g\"uzeller g\"uzeli k\"om\"ur g\"ozleriyle,\\
\c{C}akalar\i n, u\~grular\i n, e\~grilerin elinde. \\
$\; $\\
$\;\;\;\;\;\;\;\;\;\;\;\;\;\;\;\;\;\;\;\;\;\;\;\;\;\;\;\;\;\;\;\;\;\;\;\;\;\;\;\;\;\;\;\;\;\;\;\;$ \"O. Hayyam
\end{verse}
$\;$

\end{table}
\end{document}